\documentclass[12pt,twoside]{article}

\usepackage{a4}
\usepackage[centertags]{amsmath}
\usepackage{amssymb}
\usepackage{amsthm}
\usepackage{ulem}
\usepackage{isolatin1}
\usepackage{color}

\usepackage{mathtools}
\usepackage{enumitem}

\newtheorem{thm}{Theorem}
       
\newtheorem{lem}[thm]{Lemma}
\newtheorem{corollary}[thm]{Corollary}

\newtheorem{fremdersatz}{Theorem}

\newenvironment{satOrig}[1]
        {\pagebreak[2] \begin{fremdersatz} {\bf #1} \quad\sl}
        { \end{fremdersatz}}

\newenvironment{proofof}[1]
        {\pagebreak[2] \vspace{-1pt}{\bf Proof#1.}  }
        {\hfill $\blacksquare$ \vspace{2pt}}

\newenvironment{einr}{\parmod 
                      \begin{list}{}
                        {\setlength{\rightmargin}{0cm}
                         \setlength{\leftmargin}{0,75cm}
                         \setlength{\labelwidth}{0cm}
                         \setlength{\parsep}{1pt}
                         \setlength{\itemsep}{1pt}
                         \setlength{\topsep}{1pt}
                         \setlength{\partopsep}{0pt}
                         \setlength{\labelsep}{0cm}
                         \setlength{\listparindent}{0pt}
                         \setlength{\itemindent}{0pt}}
                      \item[] \ignorespaces}
                     {\unskip \end{list}}

\numberwithin{equation}{section}

\def\parmod{\parskip=2pt plus1pt minus1pt}

\newcommand{\fa}{\mathcal{F}}
\newcommand{\ha}{\mathcal{H}}

\newcommand{\ma}{\mathcal{M}}

\newcommand{\nat}{{\rm I\! N}}

\newcommand{\co}{{\mathbb C}}

\newcommand{\zet}{{\mathbb Z}}

\renewcommand{\l}{\left}
\renewcommand{\r}{\right}
\newcommand{\gl}{\left\{}
\newcommand{\gr}{\right\}}
\newcommand{\kl}{\left(}
\newcommand{\kr}{\right)}
\newcommand{\el}{\left[}
\newcommand{\er}{\right]}
\newcommand{\kj}{\overline}
\newcommand{\limn}{\lim_{n\to\infty}}

\newcommand{\In}{\subseteq}
\newcommand{\mi}{\setminus}

\newcommand{\convl}{\stackrel{\chi}{\Longrightarrow}}

\renewcommand{\rho}{\varrho}
\renewcommand{\phi}{\varphi}
\renewcommand{\epsilon}{\varepsilon}

\newcommand{\beq}{\begin{equation}}
\newcommand{\eeq}{\end{equation}}
\newcommand{\beqar}{\begin{eqnarray}}
\newcommand{\eeqar}{\end{eqnarray}}
\newcommand{\beqaro}{\begin{eqnarray*}}
\newcommand{\eeqaro}{\end{eqnarray*}}
\newcommand{\bsat}{\begin{thm}}
\newcommand{\esat}{\end{thm}}
\newcommand{\bsatorig}{\begin{satOrig}}
\newcommand{\esatorig}{\end{satOrig}}
\newcommand{\blem}{\begin{lem}}
\newcommand{\elem}{\end{lem}}
\newcommand{\bkor}{\begin{corollary}}
\newcommand{\ekor}{\end{corollary}}
\newcommand{\bbew}{\begin{proofof}}
\newcommand{\ebew}{\end{proofof}}

\def\beinr{\begin{einr}}
\def\eeinr{\end{einr}}

\begin{document}

\title{Differential inequalities and a Marty-type criterion for quasi-normality}  

\author{J\"urgen Grahl, Tomer Manket and Shahar Nevo}

\maketitle

\centerline{\bf Abstract} 

We show that the family of all holomorphic functions $f$ in a
domain $D$ satisfying 
$$\frac{|f^{(k)}|}{1+|f|}(z)\le C \qquad \mbox{ for all } z\in D$$
(where $k$ is a natural number and $C>0$) is quasi-normal.
Furthermore, we give a general counterexample to show that for
$\alpha>1$ and $k\ge2$ the condition  
$$\frac{|f^{(k)}|}{1+|f|^\alpha}(z)\le C \qquad \mbox{ for all } z\in
D$$
does not imply quasi-normality. 



\section{Introduction and statement of results}

One of the key results in the theory of normal families of meromorphic
functions is Marty's theorem \cite{Marty} which says that a family
$\fa$ of meromorphic functions in a domain $D$ in the complex plane
$\co$ is normal (in the sense of Montel) if and only if the family
$\gl f^\#\;:\;f\in\fa\gr$ of the corresponding spherical derivatives
$f^\#:=\frac{|f'|}{1+|f|^2}$ is locally uniformly bounded in $D$.

A substantial (and best possible) improvement of the direction
``$\Longleftarrow$'' in Marty's theorem is due to A.~Hinkkanen
\cite{Hinkk}: A family of meromorphic (resp.~holomorphic) functions is
already normal if the corresponding spherical derivatives are bounded
on the preimages of a set consisting of five (resp.~three) elements.
(An analogous result for normal functions was earlier proved by P.
Lappan \cite{lappan}.)

In several previous papers
\cite{BarGrahlNevo,ChenNevoPang,GrahlNevo-Spherical,GN-Marty,GN-QuasiInduced,GNP-NonExplicit, 
  LiuNevoPang} we studied the question how normality (or
quasi-normality) can be characterized in terms of the more general quantity
$$\frac{|f^{(k)}|}{1+|f|^\alpha} \qquad \mbox{ where } k\in\nat,\;
\alpha>0$$ 
rather than the spherical derivative $f^\#$. 

Before summarizing the main results from these studies we would like
to remind the reader of the definition of quasi-normality and also to
introduce some notations.  

A family $\fa$ of meromorphic functions in a domain $D\In\co$ is said
to be {\it quasi-normal} if from each sequence $\gl f_n\gr_n$ in $\fa$
one can extract a subsequence which converges locally uniformly (with
respect to the spherical metric) on $D\mi E$ where the set $E$ (which
may depend on $\gl f_n\gr_n$) has no accumulation point in $D$. If the
exceptional set $E$ can always be chosen to have at most $q$ points,
yet for some sequence there actually occur $q$ such points, then we say
that $\fa$ is {\it quasi-normal of order} $q$.  

We set $\Delta(z_0,r):=\gl z\in\co:|z-z_0|<r\gr $ for the open disk
with center $z_0\in\co$ and radius $r>0$.  By $\ha(D)$ we denote the
space of all holomorphic functions and by $\ma(D)$ the space of all
meromorphic functions in a domain $D$. We write $P_f$ and $Z_f$ for
the set of poles resp.~for the set of zeros of a meromorphic function
$f$, and we use the notation ``$f_n\convl f$ (in $D$)'' to indicate
that the sequence $\gl f_n\gr_n$ converges to $f$ locally uniformly in
$D$ (with respect to the spherical metric).

The Marty-type results known so far can be summarized as follows. 

\bsatorig{} \label{PrevResultsUpperBound}
Let $k$ be a natural number, $\alpha>0$ be a real number and $\mathcal
F$ be a family of functions meromorphic in a domain $D$. Consider the
family
$$\fa^*_{k,\alpha}:=\left\{\dfrac{|f^{(k)}|}{1+|f|^\alpha}:f\in\mathcal
  F\right\}.$$
Then the following holds. 
\begin{enumerate}
\item[(a)] \cite{LiuNevoPang, GNP-NonExplicit}
If each $f\in\fa$ has zeros only of multiplicity $\ge k$ and if
$\fa^*_{k,\alpha}$ is locally uniformly bounded in $D$, then $\fa$ is
normal.  
\item[(b)] (Y. Xu \cite{Xu}) 
Assume that there is a value $w^*\in\co$ and a constant $M<\infty$
such that for each $f\in\fa$ we have 
$|f'(z)|+\dots+|f^{(k-1)}(z)|\le M$ whenever $f(z)=w^*$ and that
there exists a set $E\subset\kj\co$ consisting of $k+4$ elements such that
for all $f\in\fa$ and all $z\in D$ we have  
\beq\label{Xu-Cond}
f(z)\in E \quad \Longrightarrow  \quad\dfrac{|f^{(k)}|}{1+|f|^{k+1}}(z)\le M.
\eeq
Then $\fa$ is normal. 

If all functions in $\fa$ are holomorphic, this also holds if one
merely assumes that $E$ has at least $3$ elements. 
\item[(c)] \cite{GN-Marty}
If $\alpha>1$ and if each $f\in\fa$ has poles only of multiplicity
$\ge\frac{k}{\alpha-1}$, then the normality of $\fa$ implies that  
$\fa^*_{k,\alpha}$ is locally uniformly bounded.

This does not hold in general for $0<\alpha\le 1$. 
\end{enumerate}
\esatorig

{\bf Remarks.}
\begin{enumerate}[label=(\arabic*)]
\item 
In (a) and (b) the assumption on the multiplicities of the zeros 
resp.~the (slightly weaker) condition on the existence of the value
$w^*$ is essential. The condition
$\frac{|f^{(k)}(z)|}{1+|f(z)|^\alpha}\le C$ itself does not imply
normality. Indeed, each polynomial of degree at most $k-1$ satisfies
this condition, but those polynomials only form a quasi-normal, but
not a normal family.  
\item 
It's worthwile to mention two special cases of Theorem
\ref{PrevResultsUpperBound} (c): 
\begin{itemize}
\item 
If $\alpha\ge k+1$ and if $\fa$ is normal, then the conclusion that
$\fa^*_{k,\alpha}$ is locally uniformly bounded holds without any
further assumptions on the multiplicities of the poles. This had been
proved already by S.Y.~Li and H.~Xie \cite{LiXie}. 
\item 
If all functions in $\fa$ are holomorphic, then for any $\alpha>1$ the
normality of $\fa$ implies that $\fa^*_{k,\alpha}$ is locally
uniformly bounded \cite[Theorem~1 (c)]{GNP-NonExplicit}. 
\end{itemize}
\end{enumerate}

In this paper we further study the differential inequality
$\frac{|f^{(k)}(z)|}{1+|f(z)|^\alpha}\le C$, but this time without any
additional assumptions on the multiplicities of the zeros of the
functions under consideration. It turns out that for $\alpha=1$ (and
hence trivially for $\alpha<1$) this differential inequality implies
quasi-normality, but that this doesn't hold for $\alpha>1$.

\bsat{} \label{mainresult}
Let $k\ge 2$  be a natural number, $C>0$  and $D\subseteq\co$ a domain. Then the family 
$$\fa_k:=\gl f\in \mathcal{H}(D) : \frac{|f^{(k)}(z)|}{1+|f(z)|}\le C \gr$$ 
is quasi-normal. 
\esat

\textbf{Remarks.}
\begin{enumerate}[label=(\arabic*)]
\item 
In Theorem \ref{mainresult} we restrict to holomorphic rather than
meromorphic functions, since if a meromorphic function $f$ has a pole at
$z_0$, then $\tfrac{|f^{(k)}(z)|}{1+|f(z)|}\le C$ is clearly violated
in a certain neighborhood of $z_0$.
\item 
The result also holds for $k=1$, and in this case we can even conclude
that $\fa$ is normal. However, this is just a trivial consequence of
Hinkkanen's extension of Marty's theorem since the condition
$\tfrac{|f'(z)|}{1+|f(z)|}\le C$ clearly implies that the derivatives
$f'$ (and hence the spherical derivatives $f^\#$) are uniformly
bounded on the preimages of five finite values. 
\item 
In Theorem \ref{mainresult}, for $k\ge 2$ the order of quasi-normality can be
arbitrarily large. This is demonstrated by the sequence of the
functions 
$$f_n(z):=n \kl e^z-e^{\zeta z}\kr$$ 
(where $\zeta:=e^{2\pi i/k}$) on the strip $D:=\gl z\in\co:
-1<{\rm Re}((1-\zeta) z)<1\gr$. Indeed, $f_n^{(k)}=f_n$, so the 
differential inequality from Theorem \ref{mainresult} trivially holds,
but every subsequence of $\gl f_n\gr_n$ is not normal exactly at the
infinitely many common zeros $z_j=\tfrac{2\pi ij}{1-\zeta}\in D$ ($j\in\zet$)
of the $f_n$, so $\gl f_n\gr_n$ is quasi-normal of infinite order. 
\item 
In the spirit of Bloch's heuristic principle, one might ask for a
corresponding result for entire functions. However, since the
exponential function (and more generally, entire solutions of the
linear differential equation $f^{(k)}=C\cdot f$) satisfy the condition
$\tfrac{|f^{(k)}(z)|}{1+|f(z)|}\le C$, there doesn't seem to be a
natural analogue for entire functions. 
\item 
For $\alpha>1$ and $k\ge 2$ the condition
$\frac{|f^{(k)}(z)|}{1+|f(z)|^\alpha}\le C$  does not imply
quasi-normality. In section \ref{sec:Counterex} we will 
construct a general counterexample for arbitrary $k\ge2$,
$\alpha>1$ and $C>0$. (For $k=2$ and $\alpha=3$ we had given such a
counterexample already in \cite{GNP-NonExplicit}.) 

In fact, it turns out that this condition doesn't even imply
$Q_\beta$-normality for any ordinal number $\beta$. (For the exact
definition of $Q_\beta$-normality we refer to
\cite{Nevo-Transfinite}.) So there isn't a chance of extending Theorem
\ref{mainresult} to the case $\alpha>1$ even if one replaces the
concept of quasi-normality by a weaker concept. 

The same counterexample also shows that Theorem \ref{mainresult}
cannot be extended in the spirit of the afore-mentioned results of
Hinkkanen and Xu (Theorem \ref{PrevResultsUpperBound} (b)). More
precisely,  a condition like 
$$f(z)\in E \quad \Longrightarrow \quad\frac{|f^{(k)}|}{1+|f|}(z)\le
C$$
where $E$ is any finite subset of $\co$ does not imply
quasi-normality (and not even $Q_\beta$-normality). This is due to the
fact that this condition is even weaker than
$\tfrac{|f^{(k)}(z)|}{1+|f(z)|^\alpha}\le C'$ for suitable $C'>0$. 
\end{enumerate}

One crucial step in our proof of Theorem \ref{mainresult} consists in
using the fact that also the reverse inequality
$\frac{|f^{(k)}(z)|}{1+|f(z)|}\ge C$ implies quasi-normality
\cite{GN-QuasiInduced}. This is one of the main results from our
studies \cite{BarGrahlNevo,ChenNevoPang, GN-QuasiInduced,LiuNevoPang}
on meromorphic functions satisfying differential inequalities of the
form $\frac{|f^{(k)}|}{1+|f|^\alpha}(z)\ge C$. These investigations
were inspired by the observation that there is a counterpart to
Marty's theorem in the following sense: A family of meromorphic
functions whose spherical derivatives are bounded away from zero has
to be normal \cite{GrahlNevo-Spherical,Steinmetz}. For the sake of
completeness, we summarize the main results from those studies.

\bsatorig{} \label{PrevResultsLowerBound}
Let $k\ge 1$ and $j\ge 0$ be integers and $C>0$, $\alpha>1$
be real numbers. Let $\fa$ be a family of meromorphic functions in
some domain $D$. 
\begin{itemize}
\item[(a)] \cite{ChenNevoPang} 
If 
$$\frac{|f^{(k)}|}{1+|f|^\alpha}(z)\ge C \qquad \mbox{ for all }
z\in D \mbox{ and all } f\in\fa,$$
then $\fa$ is normal.
\item[(b)] \cite{LiuNevoPang, GN-QuasiInduced} 
If 
$$\frac{|f^{(k)}|}{1+|f|}(z)\ge C \qquad \mbox{ for all }
z\in D \mbox{ and all } f\in\fa,$$
then $\fa$ is quasi-normal, but in general not normal. 
\item[(c)] \cite{BarGrahlNevo}
If $k>j$ and
$$\frac{|f^{(k)}|}{1+|f^{(j)}|^\alpha}(z)\ge C \qquad \mbox{ for all }
z\in D \mbox{ and all } f\in\fa,$$
then $\fa$ is quasi-normal in $D$. If all functions in $\fa$ are
holomorphic, $\fa$ is quasi-normal of order at most $j-1$. (For $j=0$
and $j=1$ this means that it is normal.) 

This does not hold for $\alpha=1$ if $j\ge1$. 
\end{itemize}
\esatorig

\section{Proof of Theorem \ref{mainresult}}

We apply induction. As mentioned above, the quasi-normality (in fact,
even normality) of $\fa_1$ follows from Hinkkanen's generalization of
Marty's theorem. 

Let some $k\ge 2$ be given and assume that it is already known that (on
arbitrary domains) each of the conditions  
$$\frac{|f^{(j)}(z)|}{1+|f(z)|}\le C \qquad \mbox{ where } j\in\gl1,\dots,k-1\gr $$
%\dots, \frac{|f^{(k-1)}(z)|}{1+|f(z)|}\le C$$
implies quasi-normality. 

Let $\gl f_n\gr_n$ be a sequence in $\fa_k$ and $z^*$ an arbitrary
point in $D$.  Suppose  to the contrary that $\gl f_n\gr_n$ is not
quasi-normal at $z^*$.

{\bf Case 1:} There is an $m\in\gl1,\dots,k-1\gr$ and a subsequence
$\gl f_{n_\ell}\gr_\ell$ such that both $\gl f_{n_\ell}^{(m)}\gr_\ell$
and $\gl \frac{f_{n_\ell}^{(m)}}{f_{n_\ell}}\gr_\ell$ are normal at $z^*$. 

Then (after turning to an appropriate subsequence which we again
denote by $\gl f_n\gr_n$ rather than $\gl f_{n_\ell}\gr_\ell$)
w.l.o.g. we may assume that in a certain disk $\Delta(z^*,r)=:U$ both
sequences $\gl f_n^{(m)}\gr_n$ and $\gl \frac{f_n^{(m)}}{f_n}\gr_n$
converge {\it uniformly} (with respect to the spherical metric) to
limit functions $H\in\ha(U)\cup\gl\infty\gr$ and
$L\in\ma(U)\cup\gl\infty\gr$, respectively.

{\bf Case 1.1:} $H$ is holomorphic.

For each $n$ we choose $p_n$ to be the $(m-1)$'th Taylor polynomial of
$f_n$ at $z^*$, i.e. $p_n$ has degree at most $m-1$ and satisfies
$p_n^{(j)}(z^*)=f_n^{(j)}(z^*)$ for $j=0,\dots,m-1$. Then $f_n$ has  
the representation 
$$f_n(z)=p_n(z)+\int_{z^*}^z\int_{z^*}^{\zeta_1}\dots\int_{z^*}^{\zeta_{m-1}}
f_n^{(m)}(\zeta_m)\,d\zeta_m\dots d\zeta_1.$$
Here for $n\to\infty$
$$\int_{z^*}^z\int_{z^*}^{\zeta_1}\dots\int_{z^*}^{\zeta_{m-1}}
f_n^{(m)}(\zeta_m)\,d\zeta_m\dots d\zeta_1
\convl
\int_{z^*}^z\int_{z^*}^{\zeta_1}\dots\int_{z^*}^{\zeta_{m-1}}H(\zeta_m)\,d\zeta_m\dots
d\zeta_1=:F(z)$$ 
where $F$ is holomorphic in $U$.  Since the family of polynomials of degree at most
$m-1$ is quasi-normal (cf. \cite[Theorem A.5 ]{Schiff}), we obtain the
quasi-normality of $\gl f_n\gr_n$ at $z^*$. 

{\bf Case 1.2:} $L(z^*)\ne\infty$. 

We choose $r_0\in (0;r)$ such that $|L(z)|\le |L(z^*)|+1$ for all $z\in\Delta(z^*,r_0)=:U_0$.

Then for all $z\in U_0$ and all $n$ large enough we have
$$\frac{|f_n^{(m)}|}{1+|f_n|}(z)\le\frac{|f_n^{(m)}|}{|f_n|}(z)\le
|L(z)|+1\le |L(z^*)|+2$$
so by the induction hypothesis we obtain the quasi-normality of $\gl
f_n\gr_n$ at $z^*$. 

{\bf Case 1.3:} $H\equiv\infty$ and $L(z^*)=\infty$. (This comprises
the cases that $L\equiv\infty$ and that $L$ is meromorphic with a pole at $z^*$.)

We choose $r_0\in (0;r)$ such that $|L(z)|\ge 3$ for all $z\in
\Delta(z^*,r_0)=:U_0$. Then for sufficiently large $n$, say for $n\ge
n_0$, and all $z\in U_0$ we have
$$\l|\frac{f_n^{(m)}}{f_n}(z)\r|\ge |L(z)|-1\ge 2 
\qquad\mbox{ and  }\qquad \l| f_n^{(m)}(z)\r|\ge 2.$$
Now fix an $n\ge n_0$ and a $z\in U_0$. If $|f_n(z)|\le 1$, we get 
$$\frac{|f_n^{(m)}|}{1+|f_n|}(z)\ge \frac{|f_n^{(m)}|}{2}(z)\ge 1.$$
If $|f_n(z)|\ge 1$, we get 
$$\frac{|f_n^{(m)}|}{1+|f_n|}(z)\ge \frac{|f_n^{(m)}|}{2|f_n|}(z)\ge
1.$$
Combining both cases, we conclude that 
$$\frac{|f_n^{(m)}|}{1+|f_n|}(z)\ge 1 \mbox{ for all } z\in U_0\mbox{
  and all } n\ge n_0.$$ 
so by Theorem \ref{PrevResultsLowerBound} (b) we obtain the
  quasi-normality of $\gl f_n\gr_n$ at $z^*$. 

{\bf Case 2:} For each $j=1,\dots,k-1$ and each subsequence $\gl
f_{n_\ell}\gr_\ell$ at least one of the sequences  $\gl
f_{n_\ell}^{(m)}\gr_\ell$ and $\gl
\frac{f_{n_\ell}^{(m)}}{f_{n_\ell}}\gr_\ell$  is not normal at $z^*$.  

Then, after turning to an appropriate subsequence which we again
denote by $\gl f_n\gr_n$, by Montel's theorem for all $j=1,\dots,k-1$
we find sequences $\gl w_{j,n}\gr_n$ such that $\limn w_{j,n}=z^*$ and 
such that for each $n$ we have $\l|f_n^{(j)}(w_{j,n})\r|\le 1$ or
$\l|\frac{f_n^{(j)}}{f_n}(w_{j,n})\r|\le 1$. Both cases can be unified
by writing 
\beq\label{InitialValueEstimate}
\l|f_n^{(j)}(w_{j,n})\r|\le 1+\l|f_n(w_{j,n})\r|
\qquad\mbox{ for all } j=1,\dots,k-1 \mbox{ and all } n.
\eeq
Furthermore, since $\gl f_n\gr_n$ is not quasi-normal, hence not
normal at $z^*$, we may also assume that there is a sequence $\gl
w_{0,n}\gr_n$ such that $\limn w_{0,n}=z^*$ and
$\l|f_n(w_{0,n})\r|\le 1$ for all $n$. 

We choose $r>0$ sufficiently small such that
$\overline{\Delta(z^\ast,r)} \subseteq D$, $2r<1$ and
$\frac{4r\cdot (1+C)}{1-2r}\le 1$. Then there exists an $n_0\in\mathbb{N}$
such that for all $n\ge n_0$ and all $j=0,\dots,k-1$ we have
$w_{j,n}\in\Delta(z^\ast,r)$. 

We use the notation 
$$ M(r,f) := \max_{|z-z^\ast|\leq r}|f(z)| \qquad \mbox{ for }
f\in\ha\kl\kj{\Delta(z^\ast,r)}\kr$$  
and obtain for all $n\ge n_0$, all $j=1,\dots,k-1$ and all $z\in\overline{\Delta(z^\ast,r)}$
\beqaro
|f_n^{(j)}(z)| 
&=&  \l|f_n^{(j)}(w_{j,n})+\int_{[w_{j,n};z]} f_n^{(j+1)}(\zeta)d\zeta\r| \\
&\le& \l|f_n^{(j)}(w_{j,n})\r|+|z-w_{j,n}|\cdot\max_{\zeta\in[w_{j,n};z]}|f_n^{(j+1)}(\zeta)|\\
&\le& 1+\l|f_n(w_{j,n})\r|+ 2r\cdot M\kl r,f_n^{(j+1)}\kr,
\eeqaro
where for the last estimate we have applied \eqref{InitialValueEstimate}. 

Since this holds for any $z\in\kj{\Delta(z^\ast,r)}$, we conclude that
for all $n\ge n_0$ and all $j=1,\dots,k-1$ 
$$M\kl r,f_n^{(j)}\kr\le 1+ M(r,f_n)+2r\cdot M\kl r,f_n^{(j+1)}\kr.$$
Similarly, in view of $\l|f_n(w_{0,n})\r|\le 1$ we also have
$$M\kl r,f_n\kr\le 1+ 2r\cdot M\kl r,f_n'\kr.$$
Induction yields
\beqaro
M\kl r,f_n\kr
&\le& 1+\sum_{j=1}^{k-1} (2r)^j\cdot(1+ M(r,f_n)) +(2r)^k\cdot M\kl r,f_n^{(k)}\kr \\
&\le& \sum_{j=0}^{k-1} (2r)^j+ \sum_{j=1}^{k-1} (2r)^j\cdot
M(r,f_n)+(2r)^k\cdot C\cdot (1+M(r,f_n))\\
&\le& C+\frac{1}{1-2r} + \frac{2r\cdot(1+C)}{1-2r} \cdot M(r,f_n)\\[8pt]
&\le& C+\frac{1}{1-2r} + \frac{1}{2} \cdot M(r,f_n)
\eeqaro
hence
$$M(r,f_n)\le 2C+\frac{2}{1-2r}$$
for all $n\ge n_0$. Thus $\gl f_n\gr_{n\ge n_0}$ is uniformly bounded in
$\Delta(z^\ast,r)$, hence normal at $z^\ast$ by Montel's theorem.

This completes the proof of Theorem~\ref{mainresult}.

\section{A general counterexample}
\label{sec:Counterex}

In this section we will show that for $\alpha>1$ and $k\ge 2$ the
differential inequality $\frac{|f^{(k)}(z)|}{1+|f(z)|^\alpha}\le C$ does not imply
quasi-normality. In \cite{GNP-NonExplicit} we had already given a
counterexample for the case $k=2$, $\alpha=3$. We generalize this
example to arbitrary $k\geq 2$, $\alpha>1$ and $C>0$.

For given $k_0\ge 2$, $C>0$ and $\alpha>1$, we
construct a sequence $\gl f_n\gr_{n}$ of holomorphic functions in
$D:=\Delta(0;2)$ such that
$\tfrac{|f_n^{(k_0)}(z)|}{1+|f_n(z)|^\alpha}\le C$ for all $z\in D$
and all $n$, but $\gl f_n\gr_{n}$ is not quasi-normal in $D$.

First, take $p,q\in\mathbb{N}$ such that
$1<\tfrac{p}{q}<\min\gl\alpha;2\gr$. 
The real function $h(x):=\dfrac{1+x^{p/q}}{1+x^\alpha}$ is
continuous in $[0,\infty)$ with $\displaystyle\lim_{x\to\infty}h(x)=0$,
hence there exists an $M>0$ such that 
\begin{equation} \label{upperBound}
\dfrac{1+x^{\frac{p}{q}}}{1+x^\alpha} \le M \qquad \mbox{ for all } x\ge 0
\end{equation}
Let $g_n(z):=z^n-1$ for $n\geq1$. The zeros of $g_n$ are the $n$-th
roots of unity $z^{(n)}_\ell=e^{2\pi i\ell/n}\ $
$(\ell=0,1,\dots,n-1)$, and they are all simple, $g_n'(z^{(n)}_\ell)\ne 0$. 
We consider the functions 
$$h_n:=g_n \cdot e^{p_n},$$ 
where the $p_n$ are polynomials yet to be determined. Then 
$$h_n'=e^{p_n}(g_n'+g_n p_n')$$ and 
\begin{equation} \label{inductionBase}
h_n''=e^{p_n}\kl 2g_n'p_n'+g_np_n'^2+g_n''+g_np_n''\kr.
\end{equation}
Our aim is to choose the $p_n$ in such a way that for $\ell=0,\dots,n-1$
\beq\label{ConditionsHn}
h_n''(z^{(n)}_\ell)=h_n^{(3)}(z^{(n)}_\ell)=\ldots=h_n^{(k_0+1)}(z^{(n)}_\ell)=0.
\eeq
We first deduce several constraints on the $p_n$ that are sufficient
for (\ref{ConditionsHn}), and then -- by an elementary result on Hermite
interpolation -- we will see that it is possible to satisfy these
constraints with polynomials $p_n$ of sufficiently large degree.  

First, in order to get $h_n''(z^{(n)}_\ell)=0$, in view of
(\ref{inductionBase}) we'll require that 
\begin{equation} \label{firstConstructionStep}
  p_n'(z^{(n)}_\ell)=-\frac{g_n''(z^{(n)}_\ell)}{2g_n'(z^{(n)}_\ell)}\
  \ \ \ (\ell=0,1,\dots,n-1). 
\end{equation}

In order to proceed we need the following lemma:

\blem\label{LemmaCounterex}
For every $k\ge 2$ we have
\[h_n^{(k)}=e^{p_n}\el kg_n'p_n^{(k-1)}+
g_n\varphi_k(p_n',\dots,p_n^{(k-1)})+
\psi_k(g_n',\dots,g_n^{(k)},p_n',\dots,p_n^{(k-2)})+
g_np_n^{(k)}\er\] 
where $\varphi_k\in\co[x_1,\dots,x_{k-1}]$,
$\psi_k\in\co[y_1,\dots,y_{k},x_1,\dots,x_{k-2}]$ are polynomials.
\elem

\bbew{}
We prove the lemma by induction on $k$. The base case $k=2$ follows
from \eqref{inductionBase} with $\varphi_2(x_1)=x_1^2$ and
$\psi_2(y_1,y_2)=y_2$. Assume that the lemma holds for some
$k\ge2$. Then differentiating gives
\beqaro%\label{Counterex-Diff}
h_n^{(k+1)} 
&= e^{p_n} \Big[ & 
\underbrace{kg_n''p_n^{(k-1)}}\ + \pmb{kg_n'p_n^{(k)}} +
\underbrace{g_n'\varphi_k(p_n',\dots,p_n^{(k-1)})} \\[5pt]  
&& +\underbracket{g_n\sum_{m=1}^{k-1}\frac{\partial\varphi_k}{\partial
    x_m}(p_n',\dots,p_n^{(k-1)})\cdot p_n^{(m+1)}}\\[5pt] 
&&+\underbrace{\sum_{m=1}^{k}\frac{\partial\psi_k}{\partial
    y_m}(g_n',\dots,g_n^{(k)},p_n',\dots,p_n^{(k-2)})\cdot
  g_n^{(m+1)}} \\[5pt] 
&& + \underbrace{\sum_{m=1}^{k-2}\frac{\partial\psi_k}{
\partial x_m} (g_n',\dots,g_n^{(k)},p_n',\dots,p_n^{(k-2)}) \cdot p_n^{(m+1)}} 
+ \pmb{g_n'p_n^{(k)}}\\[5pt] 
&& +g_np_n^{(k+1)} + \underbrace{kg_n'p_n'p_n^{(k-1)}} 
+ \underbracket{\ g_np_n'\varphi_k(p_n',\dots,p_n^{(k-1)})} \\[5pt]
&& +\underbrace{p_n'\psi_k(g_n',\dots,g_n^{(k)},p_n',\dots,p_n^{(k-2)})} 
+ \underbracket{\ g_n p_n'p_n^{(k)}}\ \Big].\\[6pt]
& = e^{p_n} \cdot \Bigl[ & 
\!\! \pmb{(k+1)g_n'p_n^{(k)}}\ +  \underbracket{\ g_n\varphi_{k+1}(p_n',\dots,p_n^{(k)})}  \\[5pt]  
&& + \underbrace{\psi_{k+1}(g_n',\dots,g_n^{(k+1)},p_n',\dots,p_n^{(k-1)}}) +  g_np_n^{(k+1)} \Bigr].
\eeqaro
where 
$$\varphi_{k+1}(x_1,\dots,x_k)
:=\sum_{m=1}^{k-1}\frac{\partial\varphi_k}{\partial x_m}(x_1,\dots,x_{k-1})\cdot x_{m+1} 
+ x_1\varphi_k(x_1,\dots,x_{k-1})+ x_1 x_k,$$
\beqaro
\psi_{k+1}(y_1,\dots,y_{k+1},x_1,\dots,x_{k-1}) 
&:=& k y_2 x_{k-1}+ y_1 \varphi_k(x_1,\dots,x_{k-1})\\[5pt]
&& +\sum_{m=1}^{k}\frac{\partial\psi_k}{\partial y_m}(y_1,\dots,y_{k},x_1,\dots,x_{k-2})
\cdot y_{m+1}  \\[5pt]
&& +\sum_{m=1}^{k-2}\frac{\partial\psi_k}{\partial x_m} (y_1,\dots,y_{k},x_1,\dots,x_{k-2}) \cdot x_{m+1}\\[5pt]
&& +k y_1 x_1 x_{k-1}+ x_1 \psi_k(y_1,\dots,y_{k},x_1,\dots,x_{k-2})
\eeqaro
are indeed polynomials of the requested form. 

Hence the lemma holds for $k+1$ as well.
\ebew

Now we inductively determine the required values of 
$p_n^{(k)}(z^{(n)}_\ell)$ for $k=2,\dots,k_0$ and
$\ell=0,\dots,n-1$. For given $k\in\gl 2,\dots,k_0\gr$, let's assume
that we already know the values of 
$p_n'(z^{(n)}_\ell),\dots,p_n^{(k-1)}(z^{(n)}_\ell)$ 
that ensure
$h_n''(z^{(n)}_\ell)=\ldots=h_n^{(k)}(z^{(n)}_\ell)=0$ 
for all admissible $\ell$. (Note that the required values of 
$p_n'(z^{(n)}_\ell)$ have been found in \eqref{firstConstructionStep}.)

In order to find the values of $p_n^{(k)}(z^{(n)}_\ell)$ (which ensure
$h_n^{(k+1)}(z^{(n)}_\ell)=0$), we apply Lemma \ref{LemmaCounterex}
with $k+1$ in place of $k$ and obtain the condition
\begin{equation} \label{inductionStep}
p_n^{(k)}(z^{(n)}_\ell) = - \frac{\psi_{k+1}(g_n',\dots,g_n^{(k+1)},p_n',\dots,p_n^{(k-1)})}{(k+1)g_n'}(z^{(n)}_\ell).
\end{equation}
(Observe that evaluating the right hand side requires only the
knowledge of values of $p_n',\dots,p_n^{(k-1)}$ that have been
previously determined.) 

It is well known (see, for example \cite[p.52]{Stoer}) that for
every $n\geq 1$ the conditions \eqref{firstConstructionStep} and 
\eqref{inductionStep} (for $k=2,\dots,k_0$) can be achieved with a
polynomial $p_n$ of degree at most $n k_0 -1$.  

In this way we obtain 
\[h_n''(z^{(n)}_\ell)=\ldots=h_n^{(k_0+1)}(z^{(n)}_\ell)=0.\]
In particular, each $z^{(n)}_\ell$ is a zero of  $h_n^{(k_0)}$ of
multiplicity $\geq 2$. 

Now, the functions $\dfrac{{h_n^{(k_0)}}^q}{h_n^p}$ are
entire: $h_n^p$ is entire, and its zeros $z^{(n)}_\ell$
$(\ell=0,1,\dots,n-1)$  have multiplicity $p$,
while ${h_n^{(k_0)}}^q$ has zeros at $z^{(n)}_\ell$ of
multiplicity at least $2q >p$. 
Thus $c_n:=\max_{z\in\overline{D}} \l|\dfrac{\kl h_n^{(k_0)}\kr^q}{h_n^p}(z)\r|<\infty$.
Define now for every $n\geq 1$ 
$$f_n:=a_n\cdot h_n,$$ 
where $a_n>0$ is a large enough constant such that both 
\begin{equation} \label{anLarge}
a_n\geq \Bigg(\dfrac{c_n\cdot M^q}{C^q}\Bigg)^\frac{1}{p-q}
\qquad\mbox{ i.e. } \qquad
\dfrac{c_n}{a_n^{p-q}}\leq \Big(\frac{C}{M}\Big)^q
\end{equation}
and $f_n\convl\infty$ on $\co\backslash\partial\Delta(0;1)$; the
latter can be achieved by choosing 
$$a_n\ge \frac{n}{\min\gl |h_n(z)|: |z|\le 1-\frac{1}{n} \mbox{ or } 1+\frac{1}{n}\le |z|\le n\gr}.$$ 

Then $\gl f_n\gr_{n}$ is not quasi-normal in $D$ (as it is
not normal at each point of $\partial\Delta(0;1)$), yet satisfies
\begin{equation*} %\label{counterExmIneq}
\frac{|f_n^{(k_0)}(z)|}{1+|f_n(z)|^\alpha}\le C
\qquad\mbox{ for all } z\in D.
\end{equation*}
Indeed, for all $z\in D$ we have
$$ \kl \frac{|f_n^{(k_0)}|}{1+|f_n|^{\frac{p}{q}}} \kr ^q(z)
\le \frac{|f_n^{(k_0)}|^q}{1+|f_n|^p}(z) 
\le\frac{|f_n^{(k_0)}|^q}{|f_n|^p}(z) 
= \frac{a_n^q \cdot \l|h_n^{(k_0)}\r|^q}{a_n^p \cdot |h_n|^p} (z)
\leq \frac{c_n}{a_n^{p-q}} \leq \kl\frac{C}{M}\kr^q $$
where the last inequality is just $\eqref{anLarge}$. Therefore, 
\[ \frac{|f_n^{(k_0)}|}{1+|f_n|^{\frac{p}{q}}}(z)
\leq \frac{C}{M} \qquad\mbox{ for all } z\in D,\]
and together with $\eqref{upperBound}$ we conclude that
\[
\frac{|f_n^{(k_0)}|}{1+|f_n|^\alpha}(z) =
\frac{|f_n^{(k_0)}|}{1+|f_n|^{\frac{p}{q}}}(z) \cdot
\frac{1+|f_n|^{\frac{p}{q}}}{1+|f_n|^\alpha}(z)
\leq \frac{C}{M} \cdot M = C
\qquad\mbox{ for all } z\in D,
\]
as desired.

{\bf Remark.} Actually, we have shown something stronger: The condition
 $\frac{|f^{(k)}(z)|}{1+|f(z)|^\alpha}\le C$ doesn't even imply
$Q_\beta$-normality for any ordinal number $\beta$ since the
constructed sequence $\gl f_n\gr_n$ and all of its subsequences are
not normal at any point of the continuum $\partial \Delta(0;1)$.

\vspace{10pt}
\parbox{85mm}{\sl J\"urgen Grahl\\
University of W\"urzburg \\
Department of Mathematics  \\     
97074 W\"urzburg\\
Germany\\
e-mail: grahl@mathematik.uni-wuerzburg.de}
\hfill\parbox{77mm}{\sl Tomer Manket \\
Bar-Ilan University\\
Department of Mathematics\\
Ramat-Gan 52900\\
Israel\\
e-mail: tomer.manket@live.biu.ac.il
}
\\[8pt]

\parbox{85mm}{\sl Shahar Nevo \\
Bar-Ilan University\\
Department of Mathematics\\
Ramat-Gan 52900\\
Israel\\
e-mail: nevosh@math.biu.ac.il}

\end{document}